\documentclass[12pt]{amsart}

\input{xy}
\xyoption{all}
\usepackage{graphicx}
\usepackage{color}
\usepackage{verbatim}

\newtheorem{theorem}{Theorem} [section]
\newtheorem{prop}[theorem]{Proposition}
\newtheorem{lemma}[theorem]{Lemma}
\newtheorem{cor}[theorem]{Corollary}
\newtheorem{conjecture}[theorem]{Conjecture}
\newtheorem{question}[theorem]{Question}

\theoremstyle{definition}

\newtheorem{example}[theorem]{Example}

\theoremstyle{remark}
\newtheorem{remark}[theorem]{Remark}

\numberwithin{equation}{section}
\numberwithin{figure}{section}

\newcommand\C{{\mathbb C}}

\newcommand\N{{\mathbb N}}

\renewcommand\P{{\mathbb P}}

\newcommand\R{{\mathbb R}}

\newcommand\Z{{\mathbb Z}}

\newcommand\D{{\mathbb D}}

\newcommand\del{\partial}

\newcommand\eps{\varepsilon}

\renewcommand\phi{\varphi}

\renewcommand\O{\mathcal{O}}

\newcommand\Aut{\operatorname{Aut}}

\newcommand\PSL{\mathrm{PSL}}

\renewcommand\Re {\operatorname{Re}}

\newcommand\Res {\operatorname{Res}} 

\newcommand\M {\mathrm{M}}

\newcommand\kbar {\overline{k}}

\newcommand\ord {\operatorname{ord}}
\newcommand\hhat {\hat{h}}

\begin{document}

\title{Bifurcations, intersections, and heights}

\author{Laura De Marco}

\email{demarco@math.northwestern.edu}

\date{\today}

\begin{abstract}
In this article, we prove the equivalence of dynamical stability, preperiodicity, and canonical height 0, for algebraic families of rational maps $f_t :\P^1(\C) \to \P^1(\C)$, parameterized by $t$ in a quasi-projective complex variety.  We use this to prove one implication in the if-and-only-if statement of \cite[Conjecture 1.10]{BD:polyPCF} on unlikely intersections in the moduli space of rational maps; we present the conjecture here in a more general form.  
\end{abstract}



\maketitle

\thispagestyle{empty}

\section{Introduction}

Let $f:  V\times\P^1(\C) \to \P^1(\C)$ be an {\bf algebraic family} of rational maps of degree $d\geq 2$.  That is, $V$ is an irreducible quasi-projective complex variety, and $f$ is a morphism so that $f_t := f(t,\cdot): \P^1\to\P^1$ has degree $d$ for all $t\in V$.  Fix a morphism $a: V \to \P^1$, which we view as a {\bf marked point} on $\P^1$.  When $V$ is a curve, we will alternatively view $f$ as a rational function defined over the function field $k = \C(V)$, with $a\in \P^1(k)$.  In this article, we study the relation between dynamical stability of the pair $(f,a)$, preperiodicity of the point $a$, and the canonical height of $a$ (defined over the field $k$).  In the final section, we present the general form of a conjecture on density of ``special points" in this setting of dynamics on $\P^1$ (which includes as a special case some known statements about points on elliptic curves); see Conjecture \ref{conj}, compare \cite[Conjecture 1.10]{BD:polyPCF}.  We finish the article with the proof of one part of the conjecture, as an application of this study of stability in algebraic families.

\subsection{Stability}
The pair $(f, a)$ is said to be {\bf stable} if the sequence of iterates
	$$\{t \mapsto f_t^n(a(t))\}_{n\geq 1}$$
forms a normal family on $V$.  (Recall that a family of holomorphic maps is normal if it is precompact in the topology of uniform convergence on compact subsets; i.e., any sequence contains a locally-uniformly-convergent subsequence.)  The pair $(f,a)$ is {\bf preperiodic} if there exist integers $m> n\geq 0$ so that $f_t^m(a(t)) = f_t^n(a(t))$ for all $t\in V$.  The pair $(f,a)$ is {\bf isotrivial} if there exists a branched cover $W\to V$ and an algebraic family of M\"obius transformations $M: W\times\P^1\to \P^1$ so that $M_t\circ f_t \circ M_t^{-1}: \P^1\to \P^1$ and $M_t(a(t))\in\P^1$ are independent of $t$.  

It is immediate from the definitions that either preperiodicity or isotriviality will imply stability.  In this article, we prove the converse:

\begin{theorem}  \label{stable}
Let $f$ be an algebraic family of rational maps of degree $d\geq 2$, and let $a$ be a marked point.  Suppose $(f,a)$ is stable.  Then either $(f,a)$ is isotrivial or it is preperiodic.
\end{theorem}

This is a generalization of \cite[Theorem 2.5]{Dujardin:Favre:critical}, which itself extends \cite[Theorem 2.2]{McMullen:families}, treating the case where $a$ is a critical point of $f$.  In the study of complex dynamics, it is well known that a holomorphic family $f: X\times\P^1\to\P^1$, for $X$ any complex manifold, is dynamically stable if and only if the pair $(f, c)$ is stable for all critical points $c$ of $f$ \cite{Lyubich:stability}, \cite{Mane:Sad:Sullivan}, \cite[Chapter 4]{McMullen:CDR}.  For non-isotrivial algebraic families $f: V\times\P^1\to\P^1$, McMullen proved that dynamical stability on all of $V$ implies that all critical points are preperiodic \cite[Lemma 2.1]{McMullen:families}.  Combining this with Thurston's Rigidity Theorem, he concluded that a non-isotrivial stable family must be a family of flexible Latt\`es maps (i.e.~covered by an endomorphism of a non-isotrivial family of elliptic curves).    

\subsection{Canonical height}
One step in the proof of Theorem \ref{stable} provides an elementary geometric proof of Baker's theorem \cite[Theorem 1.6]{Baker:functionfields} on the finiteness of rational points with small height, for the canonical height $\hat{h}_f$ associated to the function field $\C(V)$, when the variety $V$ has dimension 1. In fact, we obtain his statement under the weaker hypothesis that $f$ is not isotrivial over $k = \C(V)$ (rather than assuming $f$ is non-isotrivial over any extension of $k$); see \cite[Remark 1.7(i)]{Baker:functionfields}.  The map $f$ is {\bf isotrivial over} $k=\C(V)$ if there exists an algebraic family of M\"obius transformations $M: V\times\P^1\to \P^1$ so that $M_t\circ f_t\circ M_t^{-1}$ is independent of $t$.  

\begin{theorem}  \label{Baker}
Suppose $f$ is a rational function defined over the function field $k = \C(V)$, of degree $\geq 2$, and assume that $V$ has dimension 1.  Let $\hat{h}_f$ be the canonical height of $f$.  If $f$ is not isotrivial over $k$, then there exists $b>0$ so that the set 
	$$\{a \in \P^1(k): \hat{h}_f(a) < b\}$$
is finite.
\end{theorem}

\begin{remark} 
Theorem \ref{Baker} contains as a special case the corresponding result about rational points on an elliptic curve $E$ over $k$, equipped with the N\'eron-Tate height, generally attributed to Lang and N\'eron \cite{Lang:Neron}; see \cite[Theorem III.5.4]{Silverman:Advanced} for a proof.  It was a step in proving the Mordell-Weil theorem for function fields.  (To treat this case, we project $E$ to $\P^1$ and let $f$ be the rational function induced by multiplication-by-2 on $E$.)   In this setting, one can say more:  the set of points in $\P^1(k)$ with height $< b$ that lift to rational points in $E(k)$ is finite for {\em all} $b>0$; see \cite[Theorem B.9]{Baker:functionfields}, where this is deduced from the conclusion of Theorem \ref{Baker}.
\end{remark}

The canonical height $\hat{h}_f$ was introduced in \cite{Call:Silverman}, and it satisfies $\hat{h}_f(f(a)) = d \,\hat{h}_f(a)$ when $f$ has degree $d$.  So Theorem \ref{Baker} implies that rational points of height 0 are preperiodic unless $f$ is isotrivial over $k$.  This was proved for polynomials $f$ in \cite{Benedetto:polynomial}.   When $k$ is a number field, this was observed in \cite{Call:Silverman}, and the conclusion of Theorem \ref{Baker} holds for all $b>0$.  In the function field setting, the conclusion of Theorem \ref{Baker} {\em cannot} hold for all $b>0$, since, for example, the union of all constant points $a \in \P^1(\C)$ will form an infinite set of bounded canonical height for any $f$.   In \S\ref{examples}, we provide explicit examples of functions $f$ for which we compute the sharp bound $b$ and the total number of rational preperiodic points.  Also, note that examples do exist of rational functions $f \in k(z)$ that are isotrivial but {\em not} isotrivial over $k = \C(V)$.  A necessary condition is a nontrivial automorphism group of $f$; see Example \ref{automorphism}.

Combining Theorem \ref{stable} with Theorem \ref{Baker}, we have:

\begin{theorem} \label{higher}
Suppose $f: V\times\P^1\to \P^1$ is a non-isotrivial algebraic family of rational maps, where $V$ has dimension 1.   Let $\hat{h}_f: \P^1(\kbar) \to \R$ be a canonical height of $f$, defined over the function field $k =\C(V)$.  For each $a\in \P^1(\kbar)$, the following are equivalent:
\begin{enumerate}
\item  the pair $(f,a)$ is stable; 
\item  $\hat{h}_f(a) = 0$;
\item  $(f,a)$ is preperiodic.
\end{enumerate}
Moreover, the set 
	$$\{a\in \P^1(k): (f,a) \mbox{ is stable}\}$$
is finite.  
\end{theorem}

\subsection{Application to intersection theory}

Combining Theorem \ref{stable} with Montel's Theorem on normal families, we obtain another argument for the ``easy" implication of the Masser-Zannier theorem on anomolous torsion for elliptic curves \cite{Masser:Zannier, Masser:Zannier:2}.

\begin{prop}  \label{elliptic}
Suppose that $E$ is any non-isotrivial elliptic curve defined over a function field $k = \C(V)$, where $V$ has dimension 1, and let $P$ be a point of $E(k)$.  Then the set of $t\in V$ for which the point $P_t$ is torsion on $E_t$ is infinite. 
\end{prop}

\noindent
The harder part of the Masser-Zannier theorems is the following statement:  if two points $P$ and $Q$ in $E(k)$ are independent on $E$ and neither is torsion, meaning that they do not satisfy a relation of the form $mP + nQ = 0$ with integers $m$ and $n$ not both zero, then the set of $t \in V$ for which $P_t$ and $Q_t$ are both torsion on $E_t$ is finite.  (In \cite{DWY:Lattes}, we gave a dynamical proof of this harder implication for the Legendre family $E_t$.)  

The theorems of Masser and Zannier \cite{Masser:Zannier:2}, followed by a series of analogous results in the dynamical setting (e.g.~\cite{BD:preperiodic, Ghioca:Hsia:Tucker, GHT:rational, DWY:QPer1}),  led to the development of a general conjecture about rational maps and marked points -- addressing a question first posed by Umberto Zannier, but also encompassing a case of intrinsic dynamical interest, where the marked points are the critical points of the map.  A precise statement of this conjecture appears as Conjecture \ref{conj} in Section \ref{sec:conjecture}.  In the setting of marked critical points, Baker and I formulated the conjecture in \cite{BD:polyPCF}; an error in one of our definitions is corrected here (Remark \ref{correction}).  In this article, we use Theorem \ref{stable} to give a proof of one implication of the more general statement of Conjecture \ref{conj}.  This implication reduces to proving the following statement.  

Every algebraic family $f: V\times\P^1\to \P^1$ induces a (regular) projection $V \to \M_d$ from the parameter space to the moduli space $\M_d$ of conformal conjugacy classes of maps.  We say the family $f$ has {\bf dimension N in moduli} if the image of $V$ under this projection has dimension $N$ in $\M_d$.  Since $V$ is irreducible, $f$ has dimension 0 in moduli if and only if $f$ is isotrivial.

\begin{theorem}  \label{easy DAO}
Let $f: V\times\P^1\to\P^1$ be an algebraic family of rational maps of degree $d\geq 2$, of dimension $N>0$ in moduli.  Let $a_1, \ldots, a_k$, with $k\leq N$, be any marked points.  Then the set 
	$$S(a_1, \ldots, a_k) = \bigcap_{i = 1}^k \; \{t\in V:  \; a_i(t) \mbox{ is preperiodic for } f_t\}$$
is Zariski dense in $V$.  
\end{theorem}

\begin{remark} 
Conjecture \ref{conj} asserts that the set $S(a_1, \ldots, a_k)$ with $k >N$ will be Zariski dense if and only if at most $N$ of the points $a_1, \ldots, a_k$ are dynamically ``independent" on $V$.
\end{remark}


\subsection{Idea of the proof of Theorem \ref{stable}}
The proofs of McMullen \cite[Theorem 2.2]{McMullen:families} and Dujardin-Favre \cite[Theorem 2.5]{Dujardin:Favre:critical} use crucially that the point is critical.  The first ingredient of our proof is similar to their proofs, building upon the fact that there are only finitely many non-constant morphisms from a quasi-projective algebraic curve $V$ to $\P^1\setminus\{0,1,\infty\}$.  This leads to the proof of Theorem \ref{Baker}.  As a special case of Theorem \ref{Baker}, we have:

\begin{prop} \label{bounded}
Let $f: V\times\P^1\to \P^1$ be an algebraic family of rational maps of degree $d\geq 2$, and assume that $V$ has complex dimension 1.  Fix a marked point $a: V\to \P^1$, and define $g_n: V \to\P^1$ by $g_n(t) = f_t^n(a(t))$.  If $f$ is not isotrivial, and if the degrees of $\{g_n\}$ are bounded, then there exist integers $m >n \geq 0$ so that 
	$$f_t^m(a(t)) = f^n_t(a(t))$$
for all $t\in V$.
\end{prop}

\begin{remark}  If $f$ is isotrivial, then the degrees of $g_n(t) = f_t^n(a(t))$ are bounded if and only if the pair $(f,a)$ is isotrivial.   (See Proposition \ref{isotrivial}.)
\end{remark}

For Theorem \ref{stable}, it suffices to treat the case where $V$ has complex dimension 1, so $V$ is a finitely-punctured Riemann surface.  Then each $g_n(t) = f_t^n(a(t))$ extends to a holomorphic map on the compactification of $V$.   In light of Proposition \ref{bounded}, it remains to prove that stability on $V$ implies the degrees of $\{g_n\}$ are bounded.  

The second ingredient of the proof is a study of normality and escape rates near the punctures of $V$.  The arguments given in Section \ref{escape rate} are inspired by the methods of \cite{DWY:QPer1} and \cite{DWY:Lattes}.  

From a geometric point of view, the idea to show that the degrees of $\{g_n\}$ are bounded is as follows.  Let $X$ denote (the normalization of) a compactification of $V$.  Under iteration, one would typically expect that $\deg g_n \approx d \deg g_{n-1}$ where $d = \deg f$.  Viewing $g_n$ as a curve in $X\times\P^1$ (by identifying the function with its graph), the expected degree growth fails when the graph of $g_n$ passes through the indeterminacy points of $(t,z) \mapsto (t, f_t(z))$ and thus its image contains ``vertical components" over the punctures of $V$.  Lemma \ref{degenerate bound} shows that the multiplicity of the vertical component in the image of any curve is uniformly bounded by some integer $q$.  On the other hand, normality on $V$ implies that the graphs of $g_n$ are converging over compact subsets of $V$.  Therefore, if $\deg g_n \to \infty$, the graph of $g_n$ must be fluctuating wildly near the punctures of $V$ when $n$ is large (Lemma \ref{puncture}).  But this fluctuation is controlled by Proposition \ref{slow growth}.

\subsection{Acknowledgements.}
I would like to thank Matt Baker, Xander Faber, Dragos Ghioca, Curt McMullen, Joe Silverman, Xiaoguang Wang, and Hexi Ye for helpful comments and discussions.  Hexi Ye provided ideas for the details in Section \ref{escape rate}.  The research was supported by the National Science Foundation DMS-1517080.


\bigskip
\section{Isotriviality and Theorem \ref{Baker}}
\label{bounded degree}

Throughout this section, we assume that $V$ is an irreducible quasi-projective complex variety of dimension 1; i.e., $V$ is obtained from compact Riemann surface by removing finitely many points.  Let $f: V\times\P^1\to\P^1$ be an algebraic family of rational maps of degree $d\geq 2$.  We prove Theorem \ref{Baker}.  We also prove Proposition \ref{isotrivial}, to provide a characterization of isotriviality which is used in the proof of Theorem \ref{stable}; it is not needed for the proof of Theorem \ref{Baker}.

\subsection{Isotrivial maps}
By definition, $f$ is isotrivial if there exists a family $\{M_t\}$ of M\"obius transformations, regular over a branched cover $p: W \to V$, so that $M_t\circ f_{p(t)} \circ M_t^{-1}$ is constant in $t$.  For a marked point $a: V\to \P^1$, the pair $(f,a)$ is isotrivial if, in addition, the function $M_t(a(p(t)))$ is constant on $W$.  (Note that this is well defined, even if the family $M_t$ is not uniquely determined.)  The map $f$ (or the pair $(f,a)$) is isotrivial over $k = \C(V)$ if $M$ can be chosen to be an algebraic family that is regular on $V$.

\begin{lemma} \label{over k}
Assume that $V$ has dimension 1.  If $(f,a)$ is isotrivial, and if $\{a, f(a), f^2(a)\}$ is a set of three distinct functions on $V$, then $(f,a)$ is isotrivial over $k$.
\end{lemma}

\proof
Suppose $(f,a)$ is isotrivial.  Let $M: W\times \P^1\to \P^1$ be an algebraic family of M\"obius transformations, with branched cover $p: W \to V$,  so that $R = M_w\circ f_{p(w)}\circ M_w^{-1}: \P^1\to \P^1$ and $b = M_w(a(p(w)))$ are independent of $w\in W$.  By removing finitely many points from $V$, we may assume that $p: W\to V$ is a covering map.  

Fix a basepoint $t_0\in V$ and choose a point $w_0 \in p^{-1}(t_0)$.  The choices of basepoint determine a representation 
	$$\rho_f: \pi_1(V, t_0) \to \Aut(f_{t_0}) \subset \PSL_2(\C)$$
so that $\rho_f$ is trivial if and only if $(f,a)$ is isotrivial over $k$.  Indeed, choose any $\gamma\in \pi_1(V, t_0)$ and let $\eta: [0,1] \to W$ be a lift of $\gamma$ with $\eta(0) = w_0$.  Write $w_t$ for $\eta(t)$.  Then $f_{p(w_0)} = f_{p(w_1)} = f_{t_0}$ and isotriviality imply that $M_{w_0}f_{t_0}M_{w_0}^{-1} = M_{w_1}f_{t_0}M_{w_1}^{-1} $, so $\rho_f(\gamma) := M_{w_1}^{-1}M_{w_0}$ is an automorphism of $f_{t_0}$.  The triviality of $\rho_f$ is equivalent to the statement that $M_{w_0} = M_{w_1}$ for all such paths, so that $M$ descends to a regular map $M: V\times\P^1\to \P^1$.

Now choose the basepoint $t_0$ so that $\{a(t_0), f_{t_0}(a(t_0)), f_{t_0}^2(a(t_0))\}$ are three distinct points in $\P^1$.  Fix any $\gamma\in \pi_1(V, t_0)$.  For each $n\geq 0$, we have 
\begin{eqnarray*}
\rho_f(\gamma)(f_{t_0}^n(a(t_0))) = f_{t_0}^n(\rho_f(\gamma)(a(t_0)))
	&=& f_{t_0}^nM_{w_1}^{-1}M_{w_0}(a(p(w_0))) \\
	&=& f_{t_0}^nM_{w_1}^{-1}(b) \\
	&=& f_{t_0}^n(a(p(w_1)) = f_{t_0}^n(a(t_0)) 
\end{eqnarray*}
so the full orbit of $a(t_0)$ under $f_{t_0}$ lies in the fixed set of $\rho_f(\gamma)$.  By the assumption on $a$ we deduce that $\rho_f(\gamma)$ is the identity.   Therefore, $(f,a)$ is isotrivial over $k$.
\qed

\begin{example}  \label{automorphism}
Consider the rational function $f_1(z) = z + 1/z$ or the cubic polynomial $P_1(z) = z^3 - 3 z$.  Both of these functions have $z\mapsto -z$ as an automorphism.  Conjugating by $M_t(z) = tz$ and setting $s = t^2$, we see that the families 
	$$f_s(z) = z + \frac{1}{sz} \quad \mbox{ and } \quad P_s(z) = s z^3 - 3 z$$
are isotrivial over a degree-2 extension of $k = \C(s)$.  On the other hand, neither $f$ nor $P$ is isotrivial over $k$.  This can be seen by computing the critical points ($= \pm \sqrt{1/s}$ in both examples), and observing that the critical points are interchanged by a nontrivial loop in $V = \C\setminus\{0\}$.  \end{example}

\begin{prop}  \label{isotrivial}
Suppose $V$ has dimension 1 and that $f$ is isotrivial.  Let $a: V\to \P^1$ be any marked point.  The following are equivalent: 
\begin{enumerate} 
\item  the pair $(f,a)$ is isotrivial;
\item	  the pair $(f,a)$ is stable; and
\item  the degrees of $g_n(t) = f_t^n(a(t))$ are bounded.
\end{enumerate}
\end{prop}

\proof
Since $f$ is isotrivial, there exist a finite branched cover $p: W\to V$, an algebraic family of M\"obius transformations $M: W \times\P^1\to \P^1$, and a map $R: \P^1 \to \P^1$ so that $M_t \circ f_{p(t)} \circ M_t^{-1} = R$ for all $t\in W$.  Set $s = p(t)$.  

If the pair $(f,a)$ is isotrivial, then $b = M_t(a(s))$ is also independent of $t$, so the degrees of $g_n(t) = f_s^n(a(s)) = M_t^{-1}(R^n(b))$ are clearly bounded.   Thus (1) implies (3).   In addition, note that any sequence in the set $\{R^n(b)\}_{n\geq 1}\subset \P^1$ has a convergent subsequence.  This implies the normality of $\{M_t^{-1}(R^n(b))\}_{n\geq 1}$ on $W$ which in turn implies the normality of $\{f_t^n(a(t))\}_n$ on $V$.  Thus (1) implies (2).

Now suppose that $(f,a)$ is not isotrivial, so that $b(t) = M_t(a(s))$ is nonconstant on $W$.  We observe first that the degrees of $\{g_n\}$ must be unbounded, showing that (3) implies (1).  Indeed, since $b$ is nonconstant, it extends to a surjective map of finite degree from a compactification $\overline{W}$ to $\P^1$.  Choose any point $z_0\in \P^1$ which is non-exceptional for $R$; that is, so that the set of preimages $R^{-n}(z_0)$ is growing in cardinality as $n\to\infty$.  For any $D>0$, choose $n$ so that the size of the set $R^{-n}(z_0)$ is larger than $D$.  Then there is a set $P\subset \overline{W}$ of cardinality $|P|\geq D$ so that $b(t)$ is in $R^{-n}(z_0)$ for all $t\in P$.  Then $R^n(b(t)) = z_0$ for all $t\in P$.  Taking $D$ as large as desired, this shows that the degrees of $\{t\mapsto R^n(b(t))\}_n$ are unbounded.  This in turn implies that the degrees of $g_n(t) = M_t^{-1}(R^n(b(t)))$ are unbounded.  

Continuing to assume that $(f,a)$ is not isotrivial, we also see that $b(t) = M_t(a(s))$ has only finitely many critical points in $W$ and the image $b(W)$ omits at most finitely many points in $\P^1$.  The map $R$ has infinitely many repelling cycles in its Julia set, so there must exist a $t_0\in W$ so that $b'(t_0)\not=0$ and $b(t_0)$ is a repelling periodic point of $R$ for all $t\in P$.  It follows that the sequence of derivatives $\frac{\del\;}{\del t}R^n(b(t)) |_{t=t_0}$ is unbounded; so the sequence $\{R^n(b(t)) = M_t(f_s^n(a(s)))\}_n$ cannot be a normal family on all of $W$.   Therefore, $\{g_n\}_n$ also fails to be normal, and so we have proved that (2) implies (1).  
\qed

\subsection{Finiteness of nonconstant maps}
As in the proofs of McMullen's theorem \cite[Theorem 2.2]{McMullen:families} (and specifically his Proposition 4.3) and Dujardin-Favre's theorem \cite[Theorem 2.5]{Dujardin:Favre:critical}, we will need the following statement for our proof of Theorem \ref{Baker}.  

\begin{lemma}   \label{finite}
Let $\Lambda$ be any quasiprojective, complex algebraic curve.  There are only finitely many non-constant holomorphic maps from $\Lambda$ to the triply-punctured sphere $\P^1\setminus\{0,1,\infty\}$.  The bound depends only on the Euler characteristic $\chi(\Lambda)$.  
\end{lemma}

\proof
Any holomorphic map $h: \Lambda \to \P^1\setminus\{0,1,\infty\}$ extends to a meromorphic function on a (smooth) compactification $X$ of $\Lambda$.  From Riemann-Hurwitz, the degree of $h$ is bounded by the Euler characteristic $-\chi(\Lambda)$.  Any meromorphic function on $X$ is determined by its zeros, poles, and ones; indeed, the ratio of two functions $h_1$ and $h_2$ with the same zeros and poles must be constant on $X$, and if $h_1(x) = 1 = h_2(x)$ for some $x$, then $h_1 \equiv h_2$.  Thus, there are only finitely many combinatorial possibilities for $h$. 
\qed

\subsection{Proof of Theorem \ref{Baker}.}
Let $d = \deg f \geq 2$.  The canonical height $\hat{h}_f(a)$ computes the growth rate of the degrees of $t \mapsto f_t^n(a(t))$ as $n\to \infty$.  Precisely, each $a\in \P^1(\kbar)$ determines a meromorphic function $W\to \P^1$ on a branched cover $p_a: W\to V$, where the topological degree of $p_a$ coincides with the algebraic degree of the field extension $k(a)$ over $k$.  The canonical height is computed as 
	$$\hat{h}_f(a) = \frac{1}{\deg p_a} \; \lim_{n\to\infty} \frac{1}{d^n} \deg g_n$$
for the maps $g_n: W\to \P^1$ defined by $g_n(s) = f_{p_a(s)}^n(a(s))$.    This height function $\hat{h}_f$ is characterized by two conditions \cite[Theorem 1.1]{Call:Silverman}:
\begin{enumerate}
\item the difference $|\hat{h}_f(a) - (\deg a)/(\deg p_a)|$ is uniformly bounded on $\P^1(\kbar)$; and 
\item $\hat{h}_f(f(a)) = d \, \hat{h}_f(a)$ for all $a\in \P^1(\kbar)$.
\end{enumerate}
In particular, the degrees of $\{g_n\}$ are growing to infinity if and only if the canonical height of $a$ is positive. 

Suppose there is a sequence of rational points $a_m\in \P^1(k)$, $m\geq 1$,  so that 
	$$1 > \hat{h}_f(a_m) \to 0$$
as $m\to \infty$.   For each point $a \in\P^1(k)$, the {\bf length} of the orbit of $a$ is defined to be the cardinality of the set $\{f^n(a): n\geq 0\}$ in $\P^1(k)$.  Suppose further that only finitely many of the $a_m$ have infinite orbit, and that the finite orbit lengths are uniformly bounded.  In this case, all but finitely many of the $a_m$ satisfy a finite number of equations of the form $f^n(a) = f^\ell(a)$ with $n\not= \ell$; thus the set $\{a_m\}$ will be finite.  If this holds for any such sequence, then the theorem is proved.  

We can assume, therefore, that the orbit lengths of the $a_m$ are tending to infinity with $m$ or are equal to infinity for all $m$.  

From property (1) of the height function, there exists a degree $D$ so that $\deg a \geq D$ with $a\in \P^1(k)$  implies $\hat{h}_f(a) \geq 1$.  For each $m$, choose an integer $N_m \geq 0$ so that the orbit of $a_m$ has length greater than $N_m$ and 
	$$\deg f^i(a_m) \leq D$$
for all $i\leq N_m$.  Property (2) of the height function and the condition $\hat{h}_f(a_m)\to 0$ imply that we may take $N_m \to \infty$ as $m\to \infty$.  We will deduce that $f$ must be isotrivial over $k$.  

Suppose now that $f_t$ has at least three distinct fixed points for general $t\in V$.  Remove the finitely many parameters in $V$ where these three fixed points have collisions.  Then there exists a branched cover $p: W\to V$ and an algebraic family of M\"obius transformations $M: W\times\P^1\to \P^1$ so that 
	$$R_t = M_t \circ f_{p(t)} \circ M_t^{-1}$$
has its fixed points at $0, 1,\infty$ for all $t\in W$.

Set $b_m(t) = M_t(a_m(p(t)))$ for each $m$, so that $R^i(b_m) = M(f^i(a_m))$ for all iterates.  The uniform bound of $D$ on the degrees of $\bigcup_m \{f^i(a_m): i \leq N_m\}$ implies that there is a uniform bound of $D'$ on the degrees of $\bigcup_m\{R^i(b_m): i\leq N_m\}$.  For each point $b_m$, we define
	$$S_{0,m}(n) = \{t \in W: R_t^n(b_m(t)) = 0\}$$
	$$S_{1, m}(n) = \{t \in W: R_t^n(b_m(t)) = 1\}$$
and 
	$$S_{\infty, m}(n) = \{t \in W: R_t^n(b_m(t)) = \infty\}$$
Since $\{0,1,\infty\}$ are fixed points of $R_t$ for all $t$, we have 
	$$S_{0,m}(n) \subset S_{0,m}(n+1)$$ 
for all $m$ and all $n$; similarly for $S_{1,m}(n)$ and $S_{\infty,m}(n)$.  Let
	$$S_m = S_{0,m}(N_m) \cup S_{1,m}(N_m) \cup S_{\infty,m}(N_m).$$  
For each $m$ and each $n\leq N_m$, the iterate $R^n(b_m)$ determines a holomorphic map 
	$$W\setminus S_m \longrightarrow \P^1\setminus\{0,1,\infty\}.$$ 
By construction, the degree of $R^{N_m}(b_m)$ is bounded by $D'$, so we have $|S_m| \leq 3 D'$ for all $m$.  Therefore, there is a uniform bound $B$ (independent of $m$) on the number of non-constant maps from $W\setminus S_m$ to the triply punctured sphere $\P^1\setminus\{0,1,\infty\}$ (Lemma \ref{finite}).  In other words, for each $m$, at most $B$ of the first $N_m$ iterates of $b_m$ are nonconstant.  Therefore, since $N_m\to\infty$, there exists $m_0$ so that $b_{m_0}$ has at least $2d+2$ consecutive iterates that are constant and distinct.  

\begin{lemma}  \label{sampling}
Suppose $A$ and $B$ are rational functions of degree $d$ such that $A(x_i) = B(x_i)$ for a sequence of $2d+1$ distinct points $x_1, \ldots, x_{2d+1}$.  Then $A=B$.
\end{lemma}

\proof
By postcomposing $A$ and $B$ with a M\"obius transformation, we may assume that the values $A(x_i)$ and $B(x_i)$ are finite for each $i$.  Consider the difference $F = A-B$.  Then $F$ is a rational function of degree $\leq 2d$.  But $F$ vanishes in at least $2d+1$ distinct points, so $F\equiv 0$.  
\qed

\medskip
Set $x_i= R^i(b_{m_0})$ for the consecutive indices $i$ for which $x_i$ is constant in $t$.  Then $R_t(x_i) = x_{i+1}$ for all $t$ along a set of $2d+1$ distinct points $x_i\in \P^1(\C)$.  Applying Lemma \ref{sampling}, we conclude that the rational function $R_t$ of degree $d$ is independent of $t$.  In other words, $f$ is isotrivial.   

It remains to show that $f$ is in fact isotrivial over $k$, but this follows from Lemma \ref{over k}.  Indeed, for the point $b_{m_0}$ in the preceding paragraph, we had $x_i = R^i(b_{m_0})$ independent of $t$.  In other words, the pair $(f, f^i(a_{m_0}))$ is isotrivial.  In addition, the orbit length of $f^i(a_{m_0})$ is at least $2d+1 \geq 3$.  Lemma \ref{over k} states that the pair $(f, f^i(a_{m_0}))$ must be isotrivial over $k$, so $f$ itself is isotrivial over $k$.  

Finally, suppose that $f_t$ has only 1 or 2 fixed points, generally in $V$.  For a general parameter $t_0$, choose a forward-invariant set of fixed points and preimages, consisting of at least three distinct points.  Pass to a branched cover $W\to V$ on which these points can be marked holomorphically, excluding the finitely many points where collisions occur.  For each of these three points, we define the sets $S_{i,m}(n)$ as above.  If the $i$-th point is mapped to the $j$-th point by $f_t$, then $S_i(n) \subset S_j(n+1)$ for all $n$. The rest of the proof goes through exactly the same.  This completes the proof of Theorem \ref{Baker}.

\subsection{Two examples: computing canonical height and the number of rational preperiodic points}
\label{examples}  
Consider $Q_t(z) = z^2 +t$, the family of quadratic polynomials.  This family defines a non-isotrivial rational function $Q$ over the function field $k= \C(t)$.  There is a unique point in $\P^1(k)$ with finite orbit for $Q$, namely the point $a = \infty$.  Indeed, writing $a(t) = a_1(t)/a_2(t)$ for $a \in \P^1(k) \setminus\{\infty\}$, we can compute explicitly that 
	$$\deg(Q(a)) = \left\{  \begin{array}{ll}  
				2 \deg a  & \mbox{if } \deg a_1 > \deg a_2 \\
				2 \deg a +1 & \mbox{if } \deg a_1 \leq \deg a_2 \end{array} \right.$$ 
In both cases, the image $Q(a)$ will satisfy the hypothesis of the first case.  Inductively then, we have $\deg Q^n(a) = 2^{n-1}\deg Q(a)$.  Consequently, the largest possible $b$ in the statement of Theorem \ref{Baker} is $b = 1/2$, since the set $\{a\in \P^1(k): \hhat_Q(a) = 1/2\}$ is precisely the constant points $a\in \C$, while 
	$$\left|\{ a\in\P^1(k):  \hhat_Q(a) < 1/2\}\right| = 1.$$

As a second example, consider the family of flexible Latt\`es maps,
	$$L_t(z) = \frac{(z^2 - t)^2}{4 z (z-1)(z-t)},$$
defining a non-isotrivial $L$ over the field $k = \C(t)$.  This family is the quotient of the endomorphism $P \mapsto P+P$ on the Legendre family of elliptic curves $E_t = \{y^2 = x(x-1)(x-t)\}$.  (See also \S\ref{relation} where this is discussed more.)  For this example, we computed in \cite[Proposition 1.4]{DWY:Lattes} that there are exactly 4 rational preperiodic points, namely $\{0, 1, t, \infty\}$.  We also explicitly computed the height of any starting point $a\in \P^1(k)$ in \cite[Proposition 3.1]{DWY:Lattes}.  The constant points $a \in \C\setminus\{0,1\}$ form an infinite set of points of canonical height $1/2$, while 
	$$\left|\{ a\in\P^1(k):  \hhat_L(a) < 1/2\}\right| = 4.$$
Again, $b = 1/2$ is the largest possible constant in the statement of Theorem \ref{Baker}.


\bigskip
\section{Escape rate at a degenerate parameter}
\label{escape rate}

In this section, we construct a ``good" escape-rate function associated to a pair of holomorphic maps $f: \D^*\times\P^1\to\P^1$ and $a: \D^*\to \P^1$, on the punctured unit disk $\D^* = \{t\in \C: 0 < |t| < 1\}$.  The construction follows that of \cite{DWY:QPer1} and \cite{DWY:Lattes}.  I am indebted to Hexi Ye for his assistance in the proof of Proposition \ref{slow growth}.

\subsection{The setting} \label{setting}
Throughout this section, we work in homogeneous coordinates on $\P^1$.  We assume we are given a family of homogeneous polynomial maps
	$$F_t: \C^2\to \C^2,$$
of degree $d\geq 2$, parameterized by $t\in \D = \{t \in \C: |t| < 1\}$, such that the coefficients of $F_t$ are holomorphic in $t$.  Each $F_t$ is given by a pair of homogeneous polynomials $(P_t, Q_t)$, and we define $\Res(F_t)$ to be the homogeneous resultant of the polynomials $P_t$ and $Q_t$.  Recall that the resultant is a polynomial function of the coefficients of $F_t$, vanishing if and only if $P_t$ and $Q_t$ share a root in $\P^1$.  See \cite[\S2.4]{Silverman:dynamics} for more information.  We assume further that $\Res(F_t) = 0$ if and only if $t=0$, and also that at least one coefficient of $F_0$ is nonzero.   

We use the norm 
	$$\|(z_1, z_2)\| = \max\{|z_1|, |z_2|\}$$
on $\C^2$.

\subsection{The escape-rate function} \label{escape}
Let $F_t$ be as in \S\ref{setting}.  Let $A: \D\to\C^2\setminus\{(0,0)\}$ be any holomorphic map.  Write $A_t$ for $A(t)$.  

For each $n\geq 0$, the iterate $F_t^n(A_t)$ is a pair of holomorphic functions in $t$; we define
	$$a_n = \ord_{t=0} F_t^n(A_t)$$
to be the minimum of the order of vanishing of the two coordinate functions at $t=0$; so $a_0 = 0$ and $a_n$ is a non-negative integer for all $n\geq 1$.   Set
	$$F_n(t) = t^{-a_n} F_t^n(A_t)$$
so that $F_n$ is a holomorphic map from $\D$ to $\C^2\setminus\{(0,0)\}$ for each $n$.  Our main goal in this section is to prove the following statement.

\begin{prop} \label{slow growth}
The functions
	$$G_n(t) = \frac{1}{d^n} \log\|F_n(t)\|$$
converge locally uniformly on the punctured disk $\D^*$ to a continuous function $G$ satisfying 
	$$G(t) = o(\log|t|)$$
as $t\to 0$.  
\end{prop}

\begin{remark}
In \cite{DWY:QPer1} and \cite{DWY:Lattes}, we used explicit expressions for $F_t$ to deduce that the function $G$ was {\em continuous} at $t=0$ for our examples.  It remains an interesting open question to determine necessary and sufficient conditions for the functions $G_n$ to converge uniformly to a continuous function $G$ on a neighborhood of $t=0$.  
\end{remark}

\subsection{Order of vanishing}
Let $F_t$ and $\Res(F_t)$ be defined as in \S\ref{setting}.   

\begin{lemma} \label{degenerate bound}
Let $q = \ord_{t=0}\Res(F_t)$.  
There are constants $0 < \alpha< 1 <\beta$ and $\delta>0$ such that
$$\alpha |t|^q\leq \frac{||F_t(z_1,z_2)||}{||(z_1,z_2)||^d}\leq\beta,$$
for all $(z_1,z_2)\in \C^2\backslash \{(0,0)\}$ and $0 < |t| < \delta$.
\end{lemma}

\proof
The statement of Lemma \ref{degenerate bound} is essentially the content of \cite[Lemma 10.1]{BRbook}, letting $k$ be the field of Laurent series in $t$, equipped with the non-archimedean valuation measuring the order of vanishing at $t=0$.  But to obtain our estimate with the Euclidean norm, we work directly with their proof.

By the homogeneity of $F_t$, it suffices to prove the estimate assuming $\|(z_1, z_2)\| = 1$ with either $z_1=1$ or $z_2=1$.  The upper bound is immediate from the presentation of $F$, with bounded coefficients on compact subsets of $\D$.  

Write $F = (F_1(x,y), F_2(x,y))$.  The resultant $\Res(F)$ is a nonzero element of the valuation ring $\O_k = \{z\in k: \ord_{t=0} z \geq 0\}$.  From basic properties of the resultant (e.g. \cite[Prop. 2.13]{Silverman:dynamics}), there exist polynomials $g_1, g_2, h_1, h_2 \in \O_k[x,y]$ so that 
\begin{equation} \label{Resultant 1}
	g_1(x,y)F_1(x,y) + g_2(x,y)F_2(x,y) = \Res(F) x^{2d-1}
\end{equation}
and
\begin{equation} \label{Resultant 2}
	h_1(x,y)F_1(x,y) + h_2(x,y)F_2(x,y) = \Res(F) y^{2d-1}.
\end{equation}
Setting $x = 1$, equation (\ref{Resultant 1}) shows that 
	$$\min\{\ord F_1(1,z_2), \ord F_2(1,z_2)\} \leq q$$
for any choice of $z_2 \in \O_k$.  In fact, taking 
	$$M =  2\max\{  \sup\{ |g_1(1,y)|:  |t|\leq 1/2, |y|\leq 1\},  \sup\{ |g_2(1,y)|:  |t|\leq 1/2, |y|\leq 1\}  \}$$
we may find $\alpha_1> 0$ and $0 < \delta_1 < 1/2$ so that 
	$$\min\{ \inf\{|F_1(1,z_2)|: |z_2|\leq 1\}, \inf\{|F_2(1,y)|: |z_2|\leq 1\} \} \geq |\Res(F)|/M \geq \alpha_1 |t|^q$$
for all $|t| < \delta_1$.  

Similarly, setting $y = 1$ in equation (\ref{Resultant 2}), we may define the analogous $\alpha_2$ and $\delta_2$ to estimate $F_1(z_1,1)$ and $F_2(z_1, 1)$ for any $|z_1|\leq 1$; the conclusion follows by setting $\alpha = \min\{\alpha_1, \alpha_2\}$ and $\delta = \min\{\delta_1, \delta_2\}$.  
\qed

\subsection{Proof of Proposition \ref{slow growth}}
It is a standard convergence argument in complex dynamics that the functions $d^{-n} \log \|F_t^n(A_t)\|$ converge locally uniformly in the region where $\Res(F_t)\not=0$, exactly as in \cite{Hubbard:Papadopol}, \cite{Fornaess:Sibony}, or \cite{Branner:Hubbard:1}.  To see that the functions $G_n$ converge locally uniformly on $\D^*$, we must look at the growth of the orders $\{a_n\}$ as $n\to \infty$.  From the definition of $a_n$, we have $a_0 = 0$ and 
\begin{equation} \label{growth of orders}
  a_{n+1}=d\, a_n + \ord_{t=0} F_t(F_n(t))
\end{equation}
for all $n$.  Hence by Lemma \ref{degenerate bound}, noting that $z = F_n(t)$ has norm bounded away from both 0 and $\infty$ as $t\to 0$, we find that
  $$ 0\leq k_{n+1} := a_{n+1}-d\cdot a_n\leq q.$$
Consequently, the sequence $a_n/d^n = \sum_{i = 1}^n k_i/d^i$ has a finite limit.  In particular, we may conclude that the sequence 
  $$G_n(t)=\frac{\log \|F_n(t)\|}{d^n}=\frac{\log \|F_t^n(A_t)\|}{d^n}-\frac{a_n}{d^n} \log |t|$$
converges locally uniformly (to $G(t)$) in the punctured unit disk. 

To show that $G(t)=o(\log|t|)$ it suffices to show that for any $\eps>0$, there is a constant $C$ and a $\delta >0$ such that 
  $$|G(t)|\leq \eps \,  |\log |t|| + C$$
for all $t$ in the disk of radius $\delta$. 

Fix a positive integer $N$, and define 
    $$b_n := \ord_{t=0} F_t^{n-N}(F_N(t))$$
for $n\geq N$, so that $b_N=0$ and 
  $$0\leq \ell_{n+1} := b_{n+1}-d\cdot b_n\leq q$$
by Lemma \ref{degenerate bound}.  In particular, we have 
	$$\frac{b_n}{d^n} = \sum_{i = N+1}^n \frac{\ell_i}{d^i} \leq \sum_{i=N+1}^\infty \frac{q}{d^i}$$
for all $n>N$.  By increasing $N$ if necessary, we can assume that 
 	$$\sum_{i=N+1}^\infty \frac{q}{d^i} < \eps.$$
Therefore (recalling the constants $0 < \alpha < 1$ and $\delta>0$ from Lemma \ref{degenerate bound}),
\begin{eqnarray*}
\frac {\log \| F_n(t)\|}{d^n}+\frac{b_n}{d^n}\log|t|&=& 
\frac{ \log\|F_t^{n-N}(F_N(t))\|}{d^n}\\
&=& \sum_{i=1}^{n-N}\left(\frac{\log\|F_t^{i}(F_N(t))\|}{d^{i+N}}-\frac{\log\|F_t^{i-1}(F_N(t))\|}{d^{i+N-1}}\right)+\frac{\log\|F_N(t))\|}{d^N}\\
&=& \sum_{i=1}^{n-N}\frac{1}{d^{i+N}}\log \frac{||F_t^i(F_N(t))||}{||F_t^{i-1}(F_N(t))||^d}+\frac{\log\|F_N(t))\|}{d^N}\\
&\geq&  \sum_{i=1}^{n-N}\frac{\log |t|^q + \log \alpha}{d^{i+N}}+\frac{\log\|F_N(t))\|}{d^N}\\
&\geq& \eps \log|t|+\frac{\log\|F_N(t))\|}{d^N} +  \sum_{i=N}^\infty \frac{\log \alpha}{d^i}.
\end{eqnarray*}
Let $C=\sup \left|\frac{\log\|F_N(t))\|}{d^N}\right|$ for $t$ in the disk of radius $1/2$. Then 
  $$\frac {\log \| F_n(t)\|}{d^n} \geq \eps \log|t|-C + \sum_{i=N}^\infty \frac{\log \alpha}{d^i} -\frac{b_n}{d^n}\log|t| \geq \eps \log|t|-C + \sum_{i=N}^\infty \frac{\log \alpha}{d^i},$$
for all $|t| < \delta$ and all $n\geq N$. 

For the reverse estimate, we have 
\begin{eqnarray*}
\frac {\log \| F_n(t)\|}{d^n}+\frac{b_n}{d^n}\log|t|&=& \sum_{i=1}^{n-N}\frac{1}{d^{i+N}}\log \frac{||F_t^i(F_N(t))||}{||F_t^{i-1}(F_N(t))||^d}+\frac{\log\|F_N(t))\|}{d^N} \\
&\leq &  \sum_{i=N}^\infty \frac{\log \beta}{d^i}  + \frac{\log\|F_N(t))\|}{d^N} 
\end{eqnarray*}
where $\beta>1$ is the constant from Lemma \ref{degenerate bound}.  With the same $C$ as above, we conclude that 
  $$\frac {\log \| F_n(t)\|}{d^n} \leq - \frac{b_n}{d^n}\log|t|  + \sum_{i=N}^\infty \frac{\log \beta}{d^i} + C \leq -\eps \log|t| + \sum_{i=N}^\infty \frac{\log \beta}{d^i} + C  $$
for all $|t| < \delta$ and all $n\geq N$.  Passing to the limit as $n\to\infty$, we conclude that $G(t) = o(\log|t|)$ for $t$ near 0.  This concludes the proof of Proposition \ref{slow growth}.

\subsection{Stability.}  Here, we gather some consequences of Proposition \ref{slow growth} that will be used in the proof of Theorem \ref{stable}.  Let $F_t$ be given as in \S\ref{setting}, so it induces a family of rational maps $f$ of degree $d$, parameterized by the punctured disk $\D^*$.  Let $a: \D\to \P^1$ be a holomorphic map with a holomorphic lift $A: \D\to \C^2\setminus\{(0,0)\}$.  We define the functions $F_n$ and $G_n$ as in \S\ref{escape}, with 
	$$G(t) = \lim_{n\to\infty} \frac{1}{d^n} \log \| F_n(t) \|.$$
Recall that the pair $(f,a)$ is stable on $\D^*$ if the sequence of holomorphic functions $\{g_n(t) := f_t^n(a(t))\}$ forms a normal family on $\D^*$.  

\begin{cor}  \label{zero}
Suppose the pair $(f,a)$ is stable on the punctured disk $\D^*$.  Then there exists a choice of holomorphic lift $A: \D \to \C^2\setminus \{(0,0)\}$ of $a$ so that $G \equiv 0$.  
\end{cor}

\proof
Stability of $(f,a)$ on $\D^*$ implies that $G$ is harmonic where $t\not=0$ for any choice of holomorphic lift $A$ of $a$ \cite[Theorem 9.1]{D:lyap}.  Indeed, take a subsequence of $\{g_n\}$ that converges uniformly on a small neighborhood $U$ in $\D^*$ to a holomorphic map $h$ into $\P^1$.  Shrinking $U$ if necessary, we may select the norm on $\C^2$ so that $\log\|s(\cdot)\|$ is harmonic on a region containing the image $h(U)$ in $\P^1$, where $s$ is any holomorphic section of $\C^2\setminus\{(0,0)\} \to \P^1$.  Then the corresponding subsequence of the harmonic functions $G_n$ are converging uniformly to a harmonic limit.  

Fix a choice of $A:\D\to \C^2\setminus\{(0,0)\}$, and construct the escape-rate function $G$.  The bound on $G$ from Proposition \ref{slow growth} implies that $G$ extends to a harmonic function on the entire disk.  Indeed, by a standard argument in complex analysis, we fix a small disk of radius $r$ and let $h$ be the unique harmonic function on this disk with $h = G$ on the boundary circle.  For each $\eps>0$, consider
	$$u_\eps(t) = G(t) - h(t) + \eps \log|t|$$
for $t$ in the punctured disk.  The function $u_\eps$ extends to an upper-semi-continuous function setting $u_\eps(0) =-\infty$, and so $u_\eps$ is subharmonic on the disk because it satisfies the sub-mean-value property.  Thus, $u_\eps \leq \eps \log r$ on the disk by the maximum principle.  Letting $\eps \to 0$, we deduce that $G \leq h$ on the punctured disk.  Applying the same reasoning to 
	$$v_\eps(t) = h(t) - G(t) + \eps \log|t|$$
we obtain the reverse inequality, that $h \leq G$, and therefore, $h = G$.  

The harmonic function $G$ can now be expressed locally as $\Re \eta$ for a holomorphic function $\eta$ on $\D$.  Now replace $A_t$ with $\tilde{A}_t = e^{-\eta(t)}A_t$.  Then 
	$$F_t^n(\tilde{A}_t) = e^{-d^n\eta(t)} F_t^n(A_t),$$ 
so the order of vanishing at $t=0$ is unchanged.  We obtain a new escape-rate function 
\begin{eqnarray*}
\tilde{G}(t) &=& \lim_{n\to\infty} \frac{1}{d^n}  \log\| t^{-a_n} F_t^n(\tilde{A}_t)\| \\
	&=&  \lim_{n\to\infty} \frac{1}{d^n}  \log\| t^{-a_n} e^{-d^n \eta(t)} F_t^n(A_t)\| \\
	&=&  \lim_{n\to\infty} \frac{1}{d^n}  \log\| e^{-d^n \eta(t)} F_n(t)\| \\
	&=&  G(t) + \log|e^{-\eta}| = G(t) - G(t) \equiv 0,
\end{eqnarray*}
completing the proof of the corollary.
\qed

\begin{lemma} \label{uniform bound}
Suppose that $G \equiv 0$.  The functions $\{F_n(t) = t^{-a_n}F_t^n(A_t)\}$ are uniformly bounded in $\C^2\setminus\{(0,0)\}$ on compact subsets of $\D^*$.  
\end{lemma}

\proof
Recall that the ``usual" escape rate of $F_t$ is defined by 
	$$G_{F_t}(z) = \lim_{n\to\infty} \frac{1}{d^n} \log\|F_t^n(z)\|$$
for $z$ in $\C^2$.  The local uniform convergence of the limit (on $\D^*\times \C^2\setminus\{(0,0)\}$) implies that $G_{F_t}(z)$ is continuous as a function of $(t,z)$; it is proper in $z\in \C^2\setminus \{(0,0)\}$, since the function satisfies 
	$$G_{F_t}(\alpha z) = G_{F_t}(z) + \log|\alpha|$$
for all $\alpha\in \C^*$ and all $(t,z)$.  So our desired result follows if we can show that, for each compact subset $K$ of $\D^*$, there exist constants $-\infty < c < C < \infty$ so that 
	$$c \leq  G_{F_t}(F_n(t)) \leq C$$
for all $n$ and all $t\in K$.  

Indeed, note that $G(t) = 0$ implies that $G_{F_t}(A_t) = \eta \log|t|$ for 
	$$\eta = \lim \frac{a_n}{d^n} = \lim_{n\to\infty} \sum_{i=1}^n \frac{k_i}{d^i}$$
as in the proof of Proposition \ref{slow growth}, with $0 \leq k_i \leq q$ for all $i$.  Therefore,
\begin{eqnarray*} 
G_{F_t}(F_n(t)) &=&	d^n G_{F_t}(A_t) - a_n \log|t| \\
			&=& (d^n\eta - a_n) \log|t| \\
			&=& \left( \sum_{i=n+1}^\infty \frac{k_i}{d^{i-n}} \right) \log|t| \\
			&\geq& \left(  \sum_{i=0}^\infty \frac{q}{d^i} \right) \log|t|.
\end{eqnarray*}
On the other hand, the sequence $a_n/d^n$ increases to $\eta$, so $(d^n\eta - a_n) \log|t| \leq 0$ for all $n$ and all $t\in \D^*$; therefore, $G_{F_t}(F_n(t)) \leq 0$ for all $t$ and all $n$.  
\qed

\bigskip
\section{Proof of Theorem \ref{stable}}
\label{proof}

Let $f$ be an algebraic family of rational maps of degree $d\geq 2$, parameterized by the irreducible quasi-projective complex variety $V$.  Let $a: V\to \P^1$ be a marked point.  In this section, we prove Theorem \ref{stable}. 

It suffices to prove the theorem when $V$ is one-dimensional.  For, if $t_0$ is any parameter in $V$ at which $a(t_0)$ is not preperiodic, taking any one-dimensional slice through $t_0$ on which $(f,a)$ is not isotrivial, we conclude that $(f,a)$ will not be stable on this slice.  Therefore the family of iterates cannot be normal on all of $V$.  

If $f$ is isotrivial, then the result follows immediately from Proposition \ref{isotrivial}.  If $f$ is not isotrivial and if the degrees of $g_n(t) := f_t^n(a(t))$ are bounded, then the conclusion follows immediately from Proposition \ref{bounded}.  

For the rest of the proof, assume that the degrees of $\{g_n\}$ are unbounded.  Suppose also that $(f,a)$ is stable and $f$ is not isotrivial.  We will derive a contradiction.  

Let $C$ denote the normalization of a compactification of $V$, so that we may view $f$ as a family defined over the punctured Riemann surface $C\setminus\{x_1, \ldots, x_n\}$.  The stability of $(f,a)$ implies that $\{g_n\}$ is normal on $V$.  As such, there exists a subsequence $\{g_{n_k}\}$ with unbounded degree that converges locally uniformly on $V$ to a holomorphic function $h: V\to \P^1$.  Note that $h$ might have finite degree or it may have essential singularities at the punctures $x_i$ of $V$.  In either case, we find:

\begin{lemma}  \label{puncture}
There exists a puncture of $V$ so that for any neighborhood $U$ of this puncture, and for any point $b\in \P^1$ (with at most one exception), the cardinality of $g_{n_k}^{-1}(b)\cap U$ (counted with multiplicities) is unbounded as $n_k \to \infty$.  
\end{lemma}

\proof
We apply the Argument Principle.  Fix any $b\in\P^1$ so that $h\not\equiv b$.  Choose coordinates on $\P^1$ so that $b= 0$ and $h\not\equiv \infty$.  Choose a small loop $\gamma_j$ around each puncture $x_j$ of $V$ on which $h$ has no zeros or poles. 

Consider the integral 
	$$N(\gamma_j) = \frac{1}{2\pi i} \int_{\gamma_j} \frac{h'}{h} \in \Z$$
computing the winding number of the loop $h\circ \gamma_j$ around the origin.   By uniform convergence of $g_{n_k} \to h$ on $\gamma_j$, and since $g_{n_k}$ is meromorphic, the number $N(\gamma_j)$ is equal to the difference between the number of zeros and number of poles of $g_{n_k}$ inside the circle for all $n_k$ sufficiently large.  But since $\deg g_{n_k}\to \infty$ and the functions converge uniformly to $h$ outside these small loops, the actual count of zeros and poles must be growing to infinity inside one of these circles.  
\qed

\bigskip
Fix a puncture $x_i$ of $V$ satisfying the condition of Lemma \ref{puncture}.  
Choose local coordinate $t$ on $C$ on a small disk $D$ around the puncture $x_i$ of $V$.  Choose coordinates on $\P^1$ so that the conclusion of Lemma \ref{puncture} holds for $b= 0$ and $b=\infty$.  Let $B$ be a small annulus in the disk $D$ of the form 
	$$B = \{r_0 < |t| < r_1\}$$
with $0 < 2r_0 < r_1 - r_0 < 1$.  Passing to a further subsequence if necessary, let 
	$$\kappa(n_k) = \min\{|g_{n_k}^{-1}(0)\cap  D_{r_0}|, |g_{n_k}^{-1}(\infty)\cap  D_{r_0}|\}$$
so that $\kappa({n_k})\to \infty$ with ${n_k}$.  	
	
As in Section \ref{escape rate}, choose a homogeneous polynomial lift $F_t$ of $f_t$ to $\C^2$, normalized so that the coefficients of $F_t$ are holomorphic in $t$ and not all $0$ at $t=0$.  By Proposition \ref{slow growth} and Corollary \ref{zero}, we may choose a holomorphic lift $A$ of $a$ with values in $\C^2\setminus\{(0,0)\}$ so that the escape-rate function $G$ satisfies $G(t) \equiv 0$.  From Lemma \ref{uniform bound}, we deduce that the sequence $\{F_n(t)\}$ is uniformly bounded -- away from $(0,0)$ and $\infty$ -- on the closed annulus $\overline{B} = \{r_0 \leq |t| \leq r_1\}$.  In other words, there exist constants $0 < c \leq C < \infty$ so that 
	$$c \leq \|F_n(t)\| \leq C$$
for all $n$ and all $t\in B$.  

Write 
	$$F_{n_k}(t) = (P_{n_k}(t)R_{n_k}(t), Q_{n_k}(t)S_{n_k}(t))$$
for holomorphic $P_{n_k}, R_{n_k}, Q_{n_k}, S_{n_k}$ where $P_{n_k}, Q_{n_k} \approx t^{\kappa({n_k})}$ on the disk $D$.  More precisely, there exist factors
$$P_{n_k}(t) = \prod_{i=1}^{\kappa({n_k})}(t-t_i) \qquad\mbox{ and } \qquad Q_{n_k}(t) = \prod_{i=1}^{\kappa({n_k})}(t-s_i)$$
for two disjoint sets of roots $\{t_i\}$ and $\{s_j\}$ contained in the small disk $D_{r_0}$.  Note that 
	$$|P_{n_k}(t)|, |Q_{n_k}(t)| \leq (2r_0)^{\kappa({n_k})}$$
for $|t| = r_0$, and 
	$$|P_{n_k}(t)|, |Q_{n_k}(t)| \geq (r_1-r_0)^{\kappa({n_k})}$$
for $|t| = r_1$.  The uniform bounds on $F_n(t)$ imply that 
	$$|R_{n_k}(t)|, |S_{n_k}(t)| \leq \frac{C}{(r_1-r_0)^{\kappa({n_k})}}$$
on the circle $|t| = r_1$ and 
	$$\max\{|R_{n_k}(t)|, |S_{n_k}(t)|\} \geq \frac{c}{(2r_0)^{\kappa({n_k})}}$$
for each $t$ on the circle $|t| = r_0$ and all ${n_k}$. But $\kappa({n_k})\to \infty$ with ${n_k}$ and $2r_0 < r_1-r_0$, so for large ${n_k}$, these estimates will violate the maximum principle applied to the holomorphic function $P_{n_k}'$ or $Q_{n_k}'$.  

The contradiction obtained shows that if $(f,a)$ is stable on $V$ with $f$ not isotrivial, then the degrees of $\{g_n\}$ must be bounded, returning us to the setting treated by Proposition \ref{bounded}.  This completes the proof of Theorem \ref{stable}.

\bigskip
\section{Density of intersections}
\label{DAO}

In this section, we prove Proposition \ref{elliptic} and Theorem \ref{easy DAO}.  

\subsection{Elliptic curves}  \label{relation}
We begin by explaining briefly the connection between Proposition \ref{elliptic} and the theme of this article.  Let $E_t$ be a family of smooth elliptic curves, parameterized by a quasiprojective algebraic curve $V$.  The equivalence relation $x\sim -x$ on $E_t$ induces a projection to $\P^1$.  Via this projection, the multiplication-by-$2$ map on $E_t$ descends to a rational function $f_t$ on $\P^1$ of degree $4$, called a {\bf Latt\`es map}.  (A formula for the resulting $f_t$ is shown for the Legendre family $E_t$ in \S\ref{examples}, defined there as $L_t$.)  The family $f_t$ is non-isotrivial if and only if the family $E_t$ is non-isotrivial.  A point $P_t\in E_t$ projects to a preperiodic point for $f_t$ if and only if $P_t$ is torsion.  

In Theorem \ref{Baker}, we also refer to the canonical height function:  by its definition, the N\'eron-Tate height on an elliptic curve is equal to $\frac12$ times the canonical height for the associated multiplication-by-2 Latt\`es map.  Therefore, height 0 on the elliptic curve coincides with height 0 for the rational function.  

\subsection{Proof of Proposition \ref{elliptic}.}
Let $E$ be a non-isotrivial elliptic curve defined over a function field $k = \C(X)$ for an irreducible complex algebraic curve $X$.  We view $E$ as a family $E_t$ of smooth elliptic curves, for all but finitely many $t\in X$; alternatively, we view $E$ as a complex surface, equipped with an elliptic fibration $E\to X$.  Fix $P\in E(k)$.  Then $P$ determines a section $P: X \to E$.  Composing this section $P$ with the degree-two quotient from each $E_t$, $t\in X$, to $\P^1$, we obtain a marked point $a_P: X \to \P^1$ and a non-isotrivial algebraic family of Latt\`es maps $f: V \times\P^1\to \P^1$ on a Zariski open subset $V\subset X$.  

From Theorem \ref{stable}, we know that the pair $(f, a_P)$ is stable if and only if the pair $(f, a_P)$ is preperiodic.  So either $a_P(t)$ is preperiodic for $f_t$ for all $t\in V$ (in which case $P$ is torsion on $E/k$), or the pair $(f, a_P)$ is not stable.

In this way, the proposition is a consequence of the following statement, which is a direct application of Montel's Theorem on normal families; see e.g.~\cite{Milnor:dynamics} for background on Montel's Theorem.  

Given a pair $(f,a)$, the {\bf stable set} $\Omega(f,a)\subset V$ is the largest open set on which $\{t\mapsto f_t^n(a(t))\}$ forms a normal family.  The {\bf bifurcation locus} $B(f,a)\subset V$ is the complement of $\Omega(f,a)$ in $V$.

\begin{prop} \label{infinite}
Suppose $f: V\times\P^1\to\P^1$ is an algebraic family of rational maps of degree $d\geq 2$.  Let $a:V\to \P^1$ be a marked point, and suppose that the bifurcation locus $B(f,a)$ is non-empty.  Then for each open $U\subset V$ intersecting $B(f,a)$, there are infinitely many $t\in U$ where $a(t)$ is preperiodic to a repelling cycle of $f_t$.  
\end{prop}

\proof
Fix an open set $U$ having nonempty intersection with $B(f,a)$.  Choose point $t_0$ in $B(f,a)\cap U$.  Choose three distinct repelling periodic points $z_1(t_0), z_2(t_0), z_3(t_0)$ of $f_{t_0}$ that are not in the forward orbit of $a(t_0)$.  Shrinking $U$ if necessary, the Implicit Function Theorem implies that these periodic points can be holomorphically parameterized by $t\in U$.  By Montel's Theorem, the failure of normality of $\{t\mapsto f_t^n(a(t))\}$ on $U$ implies there exists a parameter $t_1\in U$ and an integer $n_1>0$ so that $f_{t_1}^{n_1}(a(t_1))$ is an element of the set $\{z_1(t_1), z_2(t_1), z_3(t_1)\}$.  In particular, $a(t_1)$ is preperiodic for $f_{t_1}$.  Shrinking the neighborhood $U$, we may find infinitely many such parameters.  
\qed

\begin{remark} In \cite{DWY:Lattes}, we studied the {\em distribution} of the parameters $t\in X$ for which a marked point $P_t\in E_t$ is torsion, with $E_t$ the Legendre family of elliptic curves.  The set of such parameters is dense in the parameter space (in the usual analytic topology).  
\end{remark}

\subsection{Proof of Theorem \ref{easy DAO}.}
By hypothesis, the family $f: V\times\P^1\to \P^1$ has dimension $N$ in moduli; this means that the image of the induced projection $V\to \M_d$ to the moduli space of rational maps has dimension $N$.

To prove Zariski density, we need to show that for any algebraic subvariety $Y\subset V$ (possibly reducible), the complement $\Lambda = V\setminus Y$ contains a parameter $t$ at which all points $a_1(t), \ldots, a_k(t)$ are preperiodic.  Note that $\Lambda$ is itself an irreducible, quasi-projective complex algebraic variety, so that $f: \Lambda\times\P^1\to \P^1$ is again an algebraic family of rational maps, of dimension $N$ in moduli.  

Consider the marked point $a_1$.  Since $f$ projects to an $N$-dimensional family in the moduli space, with $N>0$, it follows that $(f, a_1)$ is not isotrivial on $\Lambda$.  By Theorem \ref{stable}, the pair $(f, a_1)$ is either preperiodic or it fails to be stable on $\Lambda$.  If $(f, a_1)$ is preperiodic, we set $\Lambda_1 = \Lambda$.  If $(f,a_1)$ is unstable, then Proposition \ref{infinite} shows there exists a parameter $t_1\in \Lambda$ where $a_1(t_1)$ is preperiodic to a repelling cycle of $f_{t_1}$.  If $a_1$ satisfies the equation $f^{n_1}(a_1) = f^{m_1}(a_1)$ at the parameter $t_1$, we define $\Lambda_1 \subset \Lambda$ to be an irreducible component of the subvariety defined by the equation $f_t^{n_1}(a_1(t)) = f^{m_1}_t(a_1(t))$ that contains $t_1$.  Then $\Lambda_1$ is a nonempty quasiprojective variety, of codimension 1 in $\Lambda$.  Furthermore, since the cycle persists under perturbation, the condition defining $\Lambda_1$ will also cut out a codimension 1 subvariety in the moduli space.  In other words, $\Lambda_1$ must project to a family of dimension $N-1$ in the moduli space.  By construction, $(f,a_1)$ is preperiodic on $\Lambda_1$.  

We continue inductively.  Fix $1 \leq i < k$.  Suppose $\Lambda_i$ is a quasiprojective subvariety of dimension $\geq N-i$ in moduli on which $(f, a_1), \ldots, (f, a_i)$ are preperiodic.  Since $N-i > N-k \geq 0$, the pair $(f, a_{i+1})$ is not isotrivial on $\Lambda_i$.   As above, we combine Theorem \ref{stable} with Proposition \ref{infinite} to find a parameter $t_{i+1}\in \Lambda_i$ where $a_{i+1}$ is preperiodic.  We define $\Lambda_{i+1}\subset \Lambda_i$ so that $(f, a_{i+1})$ is preperiodic on $\Lambda_{i+1}$, and the family $f:\Lambda_{i+1}\times \P^1\to \P^1$ has dimension at least $N-i-1$ in moduli.  In conclusion, all of the points $(f, a_1), \ldots, (f, a_k)$ are preperiodic on $\Lambda_k$, and $\Lambda_k$ has dimension at least $N-k \geq 0$ in moduli.  In particular, $\Lambda_k$ is non-empty, and the theorem is proved.

\bigskip
\section{A conjecture on intersections and dynamical relations}
\label{sec:conjecture}

We conclude this article with a revised statement of the conjecture from \cite{BD:polyPCF} on ``unlikely intersections" and density of ``special points" -- and we provide the proof of one implication, as an application of Theorem \ref{stable}.  Specifically, we look at algebraic families $f: V\times\P^1\to \P^1$ of dimension $N>0$ in moduli.  We prove that if an $(N+1)$-tuple of marked points is dynamically related, then the set of parameters $t\in V$ where they are simultaneously preperiodic is Zariski dense in $V$.  We conclude the article with an explanation of how this implies one implication of \cite[Conjecture 1.10]{BD:polyPCF}.  

\subsection{Density of special points}
In \cite[Conjecture 1.10]{BD:polyPCF}, we formulated a conjecture about the arrangement of postcritically-finite maps (``special points") in the moduli space of rational maps of degree $d\geq 2$.  It was presented as a dynamical analog of the Andr\'e-Oort Conjecture in arithmetic geometry, with the aim of characterizing the ``special subvarieties" of the moduli space, meaning the algebraic families $f: V\times\P^1\to\P^1$ with a Zariski-dense subset of postcritically-finite maps.  Roughly speaking, the special subvarieties should be those that are defined by (a general notion of) critical orbit relations.  We proved special cases of the conjecture, for certain families of polynomial maps, and we sketched the proof of one implication in the general case.   

If we formulate the conjecture to handle arbitrary marked points, not only critical points, then the statement encompasses recent results about elliptic curves, as in the work of Masser and Zannier (and therefore has overlap with the Pink and Zilber conjectures); see \cite{Masser:Zannier:2} and the references therein.  Evidence towards the more general result is given by \cite[Theorem 1.3]{BD:polyPCF} and the results of \cite{Ghioca:Hsia:Tucker, GHT:rational}.   Conjecture \ref{conj} presented here is, therefore, more than just an analogy with statements in arithmetic geometry.  

Let $V$ be an irreducible, quasiprojective complex algebraic variety, and let $f: V\times\P^1\to \P^1$ be an algebraic family of rational maps of degree $d\geq 2$.  For a collection of $n$ marked points $a_1, \ldots, a_n:  V \to \P^1$, we define 
	$$S(a_1, \ldots, a_n) = \bigcap_{i = 1}^n \; \{t\in V: \; a_i(t) \mbox{ is preperiodic for } f_t\}.$$
We say the marked points $a_1, \ldots, a_n$ are {\bf coincident} along $V$ if there exists a marked point $a_i$ and a Zariski open subset $V' \subset V$ so that 
	$$S(a_1, \ldots, a_n) \cap V' = S(a_1, \ldots, a_{i-1}, a_{i+1}, \ldots, a_n) \cap V'.$$
In other words, if $\{a_1(t), \ldots, a_{i-1}(t), a_{i+1}(t), \ldots, a_n(t)\}$ are all preperiodic for $f_t$ at a parameter $t\in V'$, then the remaining point $a_i(t)$ must also be preperiodic for $f_t$.  For example, if a pair $(f, a)$ is preperiodic on $V$, then any collection of points $\{a_1, \ldots, a_n\}$ containing $a$ will be coincident.    

A stronger notion than coincidence is that of the dynamical relation, requiring an $f$-invariant algebraic relation between the points $\{a_1, \ldots, a_n\}$.  A formal definition is given below in \S\ref{relations}.  

\begin{conjecture}  \label{conj}
Let $f:V\times\P^1\to\P^1$ be an algebraic family of rational maps of degree $d\geq 2$, of dimension $N>0$ in moduli.  Let $a_0, \ldots, a_N$ be any collection of $N+1$ marked points.  The following are equivalent:
\begin{enumerate}
\item  the set $S(a_0, \ldots, a_N)$ is Zariski dense in $V$; 
\item  the points $a_0, \ldots, a_N$ are coincident along $V$; and 
\item  the points $a_0, \ldots, a_N$ are dynamically related along $V$.
\end{enumerate}
\end{conjecture}

\begin{theorem}  \label{conjthm}
We have (3) $\implies$ (2) and (2) $\implies$ (1) in Conjecture \ref{conj}.
\end{theorem}

The implication (3) $\implies$ (2) will be a formal consequence of the definitions, while (2) $\implies$ (1) is presented below as an application of Theorem \ref{stable}. The remaining challenge is to show that (1) implies (3).  We expect that (1) $\implies$ (2) should be a consequence of ``arithmetic equidistribution" as in the proofs of \cite{BD:preperiodic, BD:polyPCF, Ghioca:Hsia:Tucker, GHT:rational, DWY:Lattes} when $V$ is a curve.  

\subsection{Dynamical relations}  \label{relations}
The basic example of a dynamical relation between two marked points $a, b: V\to \P^1$ is an orbit relation:  the existence of integers $n, m$ so that 
	$$f_t^n(a(t)) = f_t^m(b(t))$$
for all $t\in V$.  To allow for complicated symmetries, we will say that $N$ marked points $a_1, \ldots, a_N$ are {\bf dynamically related} along $V$ if there exists a (possibly reducible) algebraic subvariety 	$$X \subset (\P^1)^N$$
defined over the function field $k = \C(V)$, so that three conditions are satisfied:
\begin{enumerate}
\item[(R1)]  $(a_1, \ldots, a_N) \in X$; 
\item[(R2)] (invariance) $F(X) \subset X$, where $F= (f,f, \ldots, f): (\P^1)^N\to (\P^1)^N$; and 
\item[(R3)] (nondegeneracy) 
there exists an $i\in \{1, \ldots, N\}$ and a Zariski open subset $V'\subset V$ so that the projection from the specialization $X_t$ to the $i$-th coordinate hyperplane in $(\P^1_\C)^N$ is a finite map for all $t\in V'$.
\end{enumerate}

\begin{remark}  \label{correction}
In \cite{BD:polyPCF} after stating Conjecture 1.10, we off-handedly remarked that one implication of the conjecture ``follows easily from an argument mimicking the proof of Proposition 2.6 and the following observation."  The proof of Theorem \ref{easy DAO} in this article is the argument we had in mind, mimicking \cite[Proposition 2.6]{BD:polyPCF}, but the stated ``observation" was not formulated correctly.  The inclusion of condition (R3) and the argument below in \S\ref{corrected} are an attempt to correct that error.
\end{remark}

To illustrate the dynamical relation, observe that any $N$ points $a_1, \ldots, a_N$ are dynamically related if one of the pairs $(f,a_i)$ is preperiodic.  Indeed, if $a_i$ satisfies $f_t^n(a_i(t)) = f_t^m(a_i(t))$ for all $t\in V$, for some pair of integers $n\not= m \geq 0$, then we could take $X$ to be the hypersurface  
	$$\{x \in (\P^1)^N:  f^n(x_i) = f^m(x_i)\} $$
defined over the field $k = \C(V)$.  As a nontrivial example, we look at the relation arising in the Masser-Zannier theorems \cite{Masser:Zannier:2}.  Let $f_{[k]}$ be the Latt\`es map induced from multiplication by $k \in \N$ on a non-isotrivial elliptic curve $E$ over $k =\C(V)$.  Let $a_p, a_q$ be the projections to $\P^1$ of two points $p$ and $q$ in $E(k)$.  The linear relation $n\cdot p = m\cdot q$ between points $p$ and $q$ on $E$, for integers $n$ and $m$, translates into a dynamical relation in $(\P^1)^2$ defined by 
	$$f_{[n]}(x_1) = f_{[m]}(x_2).$$ 
This relation satisfies condition (R2) for $F = (f_{[k]}, f_{[k]})$ because all Latt\`es maps descended from the same elliptic curve must commute.

Since the writing of \cite{BD:polyPCF}, we learned about the results in \cite{Medvedev:thesis} which significantly simplify the form of possible dynamical relations.  In particular, Medvedev has shown that the varieties $X$ satisfying condition (R2) should depend nontrivially on only two input variables.  In other words, the rational function $f: \P^1\to \P^1$ will be {\em disintegrated} in the sense of \cite[Definition 2.20]{Medvedev:Scanlon}; see the first theorem in the introduction of \cite{Medvedev:Scanlon}, treating the case where $f$ is a polynomial.  An affirmative answer to the following question would provide a further refinement -- and simplification -- to the notion of dynamical relation, extending the results of \cite{Medvedev:Scanlon} beyond the polynomial setting.  (The work of Medvedev and Scanlon relied on Ritt's decomposition theory for polynomials \cite{Ritt:decompositions}; the analogous decomposition theory for rational functions is not completely understood.)  

\begin{question} \label{quest}
Assume that $f$ is not isotrivial, and suppose that points $a_1, \ldots, a_N$ are dynamically related.  Does there always exist a pair of indices $i, j$ (allowing possibly $i=j$) so that the point $(a_1, \ldots, a_N)$ satisfies a relation of the form
\begin{equation} \label{simple}
	A(x_i) = B(x_j),
\end{equation}
where $A, B \in k(z)$ are nonconstant rational functions that commute with an iterate of $f$?
\end{question}

\noindent
Hypersurfaces in $(\P^1)^N$ defined by relations of the form (\ref{simple}) satisfy condition (R2) in the definition of the dynamical relation because $f$ commutes with $A$ and $B$; they satisfy condition (R3) taking either coordinate $i$ or $j$.  

\subsection{Proof of Theorem \ref{conjthm}}  \label{corrected}
We begin by proving (3) $\implies$ (2); namely, that a dynamical relation among the points $a_0, \ldots, a_N$ implies that the points are coincident.  This follows from the definition of dynamical relation, and it does not depend on the number of points.

\begin{lemma} \label{relationcoincident}
Let $f: V\times\P^1\to \P^1$ be any algebraic family of rational maps.  Suppose marked points $a_0, a_1, \ldots, a_n$ are dynamically related along $V$.  Then the points $a_0, \ldots, a_n$ are coincident along $V$. 
\end{lemma}

\proof
Let $X$ denote the $(f, \ldots, f)$-invariant subvariety in $(\P^1)^{n+1}$ for the point $(a_1, \ldots, a_n)$, given in the definition of the dynamical relation.  Suppose the points are labelled so that $x_0$ is the coordinate satisfying condition (R3).  Then the projection from $X_t$ to $(\P^1_\C)^n$, forgetting the $0$-th coordinate, is finite, for all $t$ in the Zariski-open subset $V' \subset V$.  

Now let $t_0$ be any parameter in $V'$ at which $a_1(t_0), \ldots, a_n(t_0)$ are preperiodic for $f_{t_0}$.  The point $a_0(t_0)$ must lie in the fiber of $X_{t_0} \to (\P^1)^N$ over $(a_1(t_0), \ldots, a_n(t_0))$.  Invariance of $X$ implies the invariance of $X_{t_0}$, so that $f_{t_0}^m(a_0(t_0))$ lies in the fiber over $(f^m(a_1(t_0)),  \ldots, f^m(a_N(t_0)))$ for all $m \geq 1$.  The preperiodicity of the points guarantees that there are only finitely many points in the base in the orbit of $(a_1(t_0), \ldots, a_n(t_0))$, so the orbit of $a_0(t_0)$ must be contained in a finite set.  In other words, $a_0(t_0)$ is preperiodic.  This completes the proof. 
\qed

\medskip
Now assume (2), that the given points $a_0, \ldots, a_N$ are coincident.  Assume the points are labelled so that $a_0$ is the dependent point, in the sense that 
	$$S(a_0, \ldots, a_N) \cap V' = S(a_1, \ldots, a_N)\cap V'$$
for some Zariski-open subset $V'\subset V$.  Since $V$ has dimension $N$ in moduli, Theorem \ref{easy DAO} tells us that the set $S(a_1, \ldots, a_N)$ is Zariski dense in $V$.  Therefore, so is $S(a_0, \ldots, a_N)$, and the implication (2) $\implies$ (1) is proved.  

\subsection{Proof of one implication of  \cite{BD:polyPCF} Conjecture 1.10}
Let $f: V\times \P^1 \to \P^1$ be an algebraic family of rational maps of degree $d\geq 2$ and dimension $N>0$ in moduli.  We assume that all $2d-2$ critical points of $f$ are marked.  The conjecture of \cite{BD:polyPCF} states:  {\em the map $f_t$ is postcritically finite for a Zariski dense set of $t\in V$ if and only if there are at most $N$ dynamically independent critical points.}  

Let $c_1, \ldots, c_{2d-2}$ denote the marked critical points.  Assume that $f$ has at most $N$ dynamically independent critical points; in other words, given any $n>N$ marked critical points $c_{i_1}, \ldots, c_{i_n}$, there is a dynamical relation among them.  

Note that $S(c_1, \ldots, c_{2d-2})$ is precisely the set of parameters $t$ for which $f_t$ is postcritically finite.  Applying Lemma \ref{relationcoincident} repeatedly, and reordering the points as needed, there exists a Zariski-open subset $V' \subset V$ so that 
	$$S(c_1, \ldots, c_{2d-2}) \cap V' = S(c_1, \ldots, c_{2d-3}) \cap V' = \cdots = S(c_1, \ldots, c_N)\cap V'.$$
From Theorem \ref{easy DAO}, we know that $S(c_1, \ldots, c_N)$ is Zariski dense in $V$.  This proves that the postcritically-finite maps form a Zariski dense subset of $V$.

\bigskip \bigskip

\def\cprime{$'$}

\bigskip\bigskip

\end{document}